\documentclass[a4paper,12pt,legno]{article}

\usepackage[T1]{fontenc}
\usepackage{amsmath}
\usepackage{amsfonts}
\usepackage{amsthm}
\usepackage{mathrsfs}
\usepackage{amssymb}

\newtheorem{thm}{Theorem}
\newtheorem{prop}[]{Proposition}
\newtheorem{lem}[]{Lemma}

\newtheorem{cor}{Corollary}

\newcommand{\Rset}{\mathbb{R}}
\newcommand{\Cset}{\mathbb{C}}

\begin{document}

\title{Remarks on restricted Nevanlinna transforms\footnote{Research partially funded by the University of Wroc³aw
 grant no \textbf{2242/W/IM/09.}}}

\author{Lech Jankowski and Zbigniew J. Jurek\footnote{Corresponding
Author.}}

\date{Appeared in: \emph{Demonstratio Math.} \textbf{45} no 2, 2012, pp. 297-307.}
\maketitle

\begin{quote} \textbf{Abstract.}
The Nevalinna transform $K_{a,\rho}(z)$ of a positive measure $\rho$
and a constant $a$, plays an important role in complex analysis and
-- more recently -- in the context of the boolean convolution. We
show here that its restriction to the imaginary axis,
$k_{a,\,\rho}(it)$, can be expressed as the Laplace transform of the
Fourier transform (a characteristic function) of $\rho$.
Consequently, $k_{a\,\rho}$ is sufficient for the unique
identification of the measure $\rho$ and the constant $a$. Finally,
we identify a relation between the free additive Voiculescu
$\boxplus$ and boolean $\uplus$ convolutions.



\emph{Abbreviated title:} Nevanlinna transform.
\end{quote}

\medskip
\medskip
\newpage
The Cauchy $G(z)$ and the Nevanlinna $K(z)$ transforms play an important role
in complex analysis and free probability. They are given as follows:
\[
G_m\,(z):=\int_{\Rset}\,\frac{1}{z-x}\,m(dx), \ \ K_{a,\, \rho}
(z):= a + \int_{\Rset}\,\frac{1+zx}{z-x}\,\rho(dx), \ \ z \in
\Cset\setminus \Rset, \ \  (\ast)
\]
for some finite measures $m$ and $\rho$ and constants $a$. In order
to retrieve the measure $m$ from $G_m$ one uses the classical
inversion formula
\begin{equation*}
m([a,b])= -\lim_{y\to 0}\frac{1}{\pi}\int_a^b \Im \,G_m(x+iy)dx, \
\mbox{provided} \ \ m(\{a,b\})=0;
\end{equation*}
cf. Akhiezer (1965), p. 125 or Lang (1975), p. 380, Bondesson
(1992). Thus, $G_m$ uniquely determines $m$. It is important to
stress that the above inversion requires one to know the Cauchy
transform in strips $\{x+iy : x \in \Rset,\, 0<y<\epsilon \}$ for
some $\epsilon >0$. Jurek (2006) demonstrates that the values of
$G_m(it), t\neq 0$, are sufficient to identify $m$, using a simple
argument of exponentiation of measures; also cf. Proposition 1
below. Of course, as holomorphic functions  $G_m$ and $K_{a\,,
\rho}$ are determined by their values on sets having a condensation
point, but the proof in Jurek (2006) is notable for avoiding the use
of structural theorems from complex analysis.

This paper is an application of the general idea (conjecture) that
many transforms in complex analysis and, in particular, in the area
of the free probability, are some functionals of the standard
Laplace and Fourier transforms when suitably restricted to the
imaginary line.

In particular, we will show that the measure $\rho$, in the
Nevalinna transform, can be retrieved from values $K_{a,\,\rho}(it),
t\neq 0$, using the classical (standard) Fourier and Laplace
transforms, after restricting $K_{a,\,\rho}$ to the imaginary axis
without the origin; cf. Theorem 1 (The inversion formula). Then we
illustrate the inversion formula by an example. Finally we derive a
relation between the so-called \emph{boolean convolution} $\uplus$,
introduced by Speicher and Woroudi (1997), and the Voiculescu
convolution $\boxplus$ (Proposition 2); cf. Acknowledgement below.
Finally, Remark 2 identifies a challenging open problem.

\textbf{1. Notations, results and an example.} For a real constant
$a$ and a finite Borel measure $\rho$ on the real line, \emph{the
restricted Nevanlinna transform} is defined by
\begin{equation}
k_{a,\,\rho}(it):= a + \int_{\Rset}\,\frac{1+itx}{it-x}\,\rho(dx), \
\ \ \mbox{for}\ \ t \neq 0,
\end{equation}
and similarly, \emph{the restricted Cauchy transform}, by
\begin{equation}
g_{\rho}(it):= \int_{\Rset}\,\frac{1}{it-x}\,\rho(dx), \ \ \
\mbox{for}\ \ t \neq 0;
\end{equation}
comp. the equation $(\ast)$ above. Let us recall also that the
Fourier transform (the characteristic function) $\hat{\mu}$ of a
measure $\mu$ is given by
\begin{equation}
\hat{\mu}(t):=\int_{\Rset}\,e^{itx}\mu(dx), \ \ t\in \Rset.
\end{equation}
and the Laplace transform of a function $h: (0,\infty)\to \Cset$, or
of a measure $m$ is given by
\begin{equation}
\mathfrak{L}[h; \lambda]:=\int_0^{\infty}\,h(x)e^{-\lambda\,x}\,dx ,
\ \ \mathfrak{L}[m\,;
\lambda]:=\int_0^{\infty}\,e^{-\lambda\,x}\,m(dx), \ \ \lambda>0
\end{equation}
where $\lambda$ is a such that those integral exist; cf. Gradshteyn
and Ryzhik (1994), Chapter 17, for examples of those transforms and
their inverses.

We begin by stating our main result, showing how to obtain the measures $\rho$ knowing
only their restricted Nevanlinna transforms.
Below, $\Re z, \ \Im z,\  \overline{z}$ denote the real part, the imaginary part and
the conjugate of a complex $z\in\Cset$,
respectively.

\begin{thm}
(The inversion formula.) \ For the restricted Nevalinna transform
$k_{a,\,\rho}$  we have that: $a =\Re k_{a,\ \rho}(i), \ \
\rho(\Rset)= -\Im k_{a,\ \rho}(i) $; and  the identity
\begin{equation*}
\mathfrak{L}[\hat{\rho}; \,\,w ]=
\int_{0}^{\infty}\hat{\rho}(r)e^{-w\,r}dr=
\frac{ik_{a,\,\rho}(-iw)-i \Re k_{a,\, \rho}(i)-w \Im
k_{a,\,\rho}(i)}{w^{2}-1} \,
\end{equation*}
holds for $w>0$. In particular, the constant $a$ and the measure
$\rho$ are uniquely determined by the functional $k_{a,\,\rho}$ .
\end{thm}
Since part of the above right-hand side formula can be viewed as
Laplace transform of some exponential functions we get
\begin{cor}
For the restricted  Nevanlinna  functional $k_{a,\,\rho}$ and $w>1$
we have
\begin{equation*}
\int_{0}^{\infty}\big[\hat{\rho}(r)-\frac{1}{2}\big(i\,k_{a,\,
\rho}(i)\,e^{-x}+ \overline{i\,k_{a,\, \rho}(i)}\,e^x \big)\big]
e^{-w\,r}dr= \frac{ik_{a,\,\rho}(-iw)}{w^{2}-1} .
\end{equation*}
In particular, if $a=0$ and $ \nu$ is a probability measure then for
$k_{0,\nu}$ we get
\begin{equation*}
\int_{0}^{\infty}(\hat{\nu}(r)- \cosh r)\, e^{-w\,r}dr=
\frac{ik_{0,\,\nu}(-iw)}{w^{2}-1}, \ \ w>1\, .
\end{equation*}
\end{cor}
\begin{prop} For a finite measure $\rho$ and its restricted Cauchy
transform $g_{\rho}$ we have
\begin{equation*}
\mathfrak{L}[\hat{\rho}; \,\,w ]=\overline{i\,g_{\rho}(iw)}, \ w\neq
0,
\end{equation*}
that is, to retrieve $\rho$ one needs to invert Laplace transform of
$\hat{\rho}$ and then invert the Fourier transform.
\end{prop}

\medskip
Hence we conclude that the values of restricted Cauchy transform
$g_{\rho}(iw), \ w \neq 0,$ uniquely determine the measure $\rho$.
That fact was already established in Jurek (2006) but not explicitly
as it is in the above Proposition 1.

In the following example we show explicitly that shifted reciprocals
of restricted Cauchy transforms of discrete measures correspond to
restricted Nevannlina transforms; see the formula (5) below.

\textbf{Example.} For a set $\textbf{b}=\{b_1, b_2, ..., b_m\}$  of
distinct real numbers let us define a discrete probability measure
$\mu_{\textbf{b}}:=\frac{1}{m}\sum_{j=1}^m\delta_{b_j}$ and the
canonical polynomial $P_{\textbf{b}}(z)=\prod_{j=1}^m(z-b_j)$. If $
\{\xi_1,\xi_2,...,\xi_{m-1}\}$ is the set of zeros of the polynomial
$P'_{\textbf{b}}(z)$ (the derivative of $P$) then we have
\begin{equation}
it-\frac{1}{G_{\mu_{\textbf{b}}}(it)}= it -
\frac{m}{\sum_{j=1}^m\frac{1}{it-b_j}}= a_{\textbf{b}} +\int_{\Rset}
\frac{1+itx}{it-x} \rho_{\textbf{b}}(dx), \ \ t\neq 0,
\end{equation}
where
\begin{multline}
\alpha_k:= - m \frac{P(\xi_k)}{P''(\xi_k)}= m\,\big[
\sum_{j=1}^m\frac{1}{(\xi_k-b_j)^2}- (\sum_{j=1}^m
\frac{1}{\xi_k-b_j})^2\big]^{-1}>0 \ \ \ \ \ \ \ \
\\
a_{\textbf{b}}:=\frac{b_1+b_2+...+b_m}{m}
-\sum_{j=1}^{m-1}\frac{\alpha_j\,\xi_j}{1+\xi^2_j}, \ \
\rho_{\textbf{b}}(dx):=\sum_{j=1}^{m-1}\frac{\alpha_j}{1+\xi^2_j}\delta_{\xi_j}(dx).
\end{multline}

Note that the procedure described in the Example can be iterated.
Namely, in the second step we may start with the probability measure
concentrated on the roots $\xi_j , \ j=1,2,...,m-1$, and so on.

Recall that \emph{the self-energy functional} $E_{\mu}$ of the
probability measure  $\mu$ is defined as follows
\begin{equation}
E_{\mu}(z)= z - \frac{1}{G_{\mu}(z)} , \  \ z\in\Cset\setminus\Rset.
\end{equation}
\noindent Similarly to the above relation, we refer to $e_{\mu}(it):=E_{\mu}(it), t\neq
0$, as \emph{a restricted self-energy functional}.

To express $a$ and $\rho$ in terms of $\mu$ using only the restricted
functionals we make use of the following corollary:

\begin{cor}
For a probability measure $\mu$ let
\begin{equation}
z_{\mu}:=-g_{\mu}(i)=c_{\mu}+i\,d_{\mu}\equiv
\int_{\Rset}\frac{x}{1+x^2}\mu(dx)+i\,\int_{\Rset}\frac{1}{1+x^2}\mu(dx)\in\Cset.
\end{equation}
If $e_{\mu}(it)=k_{a,\,\rho}(it)$, for all $t\neq 0$, then the
constants $a$ and $\rho(\Rset)$ are given by formulae
\begin{equation}
a=\frac{c_{\mu}}{|z_{\mu}|^2} \ \ \  \mbox{and} \ \ \ \rho(\Rset)=
\frac{d_{\mu}}{|z_{\mu}|^2} - 1>0 \, ,
\end{equation}
and the Fourier transform $\hat{\rho}$ satisfies the equation
\begin{equation}
\mathfrak{L}[|z_{\mu}|^2\hat{\rho}(x)-\frac{1}{2}(\overline{z_{\mu}}\,\,e^x+
z_{\mu}\,e^{-x}); \,\,w ]= \frac{1}{(w^2-1)\,i\,g_{\mu}(-iw)}, \ \
w>1.
\end{equation}
\end{cor}
Since for any probability measures $\mu$ and $\nu$ there exists a unique
probability measure $\gamma$ such that
\begin{equation}
E_{\mu}(z)+E_{\nu}(z)=E_{\gamma}(z),
\end{equation}
we call it \emph{the boolean convolution} and denote it by
$\gamma=\mu\uplus\nu$; for more details cf. Speicher - Woroudi
(1997) and references therein.

\textbf{Remark 1.} Boolean convolution has the property that
\emph{all} probability measures are $\uplus$-infinitely divisible.
The \emph{max}-convolution also has that feature because \emph{for
each distribution function F, $F^{1/n}$ (the n-th root) is also a
distribution function and taking independent identically distributed
(as $F^{1/n}$ ) r.v. $X_{n, 1}, X_{n,2},...,X_{n,n}$, we see that
$\max \{X_{n,1},..., X_{n,n}\}$ has the distribution function F.}

\medskip
For a probability measure $\mu$, let
\begin{equation}
F_{\mu}(z):=\frac{1}{G_{\mu}(z)}, \ \ z \in \Cset\setminus \Rset, \
\ \mbox{and} \ \ \ V_{\mu}(z):=F_{\mu}^{-1}(z)-z, \ \ z\in
\mathcal{D} \subset \Cset,
\end{equation}
where $\mathcal{D}$ is the so called Stolz angle in which the inverse
$F_{\mu}^{-1}$ exists; cf. Bercovici-Voiculescu (1993) and
references therein. Since for any probability measures $\mu$ and $\nu$
there exists a unique probability measure $\gamma$ such that
\begin{equation}
V_{\mu}(z)+V_{\nu}(z)=V_{\gamma}(z),
\end{equation}
we call it \emph{the Voiculescu convolution} and denote it by
$\gamma=\mu\boxplus\nu$; cf. Bercovici-Voiculescu (1993) and
references therein. A relation between $\boxplus$-infinite
divisibility and some random integrals with respect to classical
L\'evy processes is given in Jurek (2007), Corollary 6.

Here are some unexpected relations between the Voiculescu $\boxplus$
and the boolean $\uplus$ operations on probability measures; cf.
Lenczewski (2007), Proposition 2.1 and the Acknowledgements below.

\begin{prop}
For probability measures $\mu_1$ and $\mu_2$ there exist unique
probability measures $\nu_1, \nu_2$ such that
\begin{equation*}
F_{\mu_1}(F_{\nu_1}(z))=F_{\mu_2}(F_{\nu_2}(z))=F_{\mu_1\boxplus\mu_2}(z),
\, \, z \in \Cset^+.
\end{equation*}
Furthermore, the above measures satisfy the equation
$\nu_1\uplus\nu_2=\mu_1\boxplus\mu_2$.
\end{prop}
\begin{cor}
For $n\ge 2$ and probability measures $\mu_1, \mu_2, ..., \mu_n$
there exist unique probability measures $\nu_1, \nu_2, ..., \nu_n$
such that $F_{\mu_1}(F_{\nu_1}(z))=F_{\mu_2}(F_{\nu_2}(z))=...=
F_{\mu_n}(F_{\nu_n}(z))=F_{\mu_1\boxplus\mu_2\boxplus...\boxplus\mu_n}(z),
\ z\in\Cset^+$. Furthermore, the above measures satisfy the equation
$(\nu_1\uplus\nu_2 \uplus...\uplus\nu_n)^{\uplus 1/(n-1)} =
\mu_1\boxplus\mu_2\boxplus...\boxplus\mu_n$.
\end{cor}
\textbf{Remark 2.} The two identities below, involving $\uplus$ and
$\boxplus$ might be of an interest in  themselves. More importantly,
finding \underline{\emph{real analytic}} proofs of them seems  to be
very challenging.

(a) For probability measures $\mu$ and $\nu$ there exists a unique
measure $\mu\uplus\nu$ such that
\begin{equation*}
\frac{1}{\int_\Rset\frac{1}{1-itx}\,\mu(dx)}\ -\ 1 \ + \
\frac{1}{\int_\Rset\frac{1}{1-itx}\,\nu(dx)} \ -\ 1
=\frac{1}{\int_\Rset\frac{1}{1-itx}\,\mu\uplus\nu(dx)}\ - \ 1, \ \
\end{equation*}
for $t\in \Rset$; cf. Theorem 2 and Remark 1.1.1 in Jurek (2006) for
other forms of the above formula and some comments.

(b) For measures $\mu_1$ and $\mu_2$ there exist unique measures
$\nu_1,\,\nu_2$ and $\mu_1\boxplus\mu_2$ such that for their
restricted Cauchy transforms we have
\begin{equation*}
g_{\nu_1}(it)\int_{\Rset}\frac{1}{1-xg_{\nu_1}(it)}\mu_1(dx)=g_{\mu_1\boxplus\mu_2}(it)=
g_{\nu_2}(it)\int_{\Rset}\frac{1}{1-xg_{\nu_2}(it)}\mu_2(dx),
\end{equation*}
for all $t\neq 0$; cf. Biane (1998), Chistyakov and Goetze (2005).
(Using Proposition 1 we may express the above identity in terms of
classical Laplace and Fourier transforms.)

\medskip
\textbf{2. Auxiliary results and proofs.} Note that
\begin{equation*}
\overline{g_m(it)}=g_m(-it), \ \ \overline{k_{\rho}(it)}=k_{a,\
\rho}(-it), \ \  \overline{e_{\mu}(it)}=e_{\mu}(-it), \ t \neq 0,
\end{equation*}
which allows us to consider those function only on the positive half-line.

\emph{Proof of Theorem 1.} From (1) we get
\begin{equation}
k_{a,\ \rho}(i)= a -i\rho(\Rset).
\end{equation}
Further, since
\[
\frac{1+itx}{it-x}=\frac{1-t^2}{it-x}-it
\]
we  infer from (1) and (2) that
\begin{equation}
k_{a,\,\rho}(it)=a+(1-t^{2})g_{\rho}(it)-it\rho(\mathbb{R}), \ \
g_{\rho}(it)=\frac{k_{a,\,\rho}(it)-a+it\rho(\Rset)}{1-t^2}\,.
\end{equation}
On the other hand, in Jurek (2006) on p. 189, it was noticed that
\begin{equation*}
\int_{0}^{\infty}\hat{\rho}(ts)e^{-s}ds =
\frac{1}{it}g_{\rho}(\frac{1}{it}), \ \ t\neq 0, \ \ \mbox{and}  \ \
\lim_{t\to 0}\frac{1}{it}g_{\rho}(\frac{1}{it})=\rho(\Rset).
\end{equation*}
This property, along with (14) and (15), yields
\begin{multline*}
\int_{0}^{\infty}\hat{\rho}(ts)e^{-s}ds =
\frac{1}{it}\,\frac{(k_{a,\,\rho}(\frac{1}{it})-\Re
k_{a,\,\rho}(i)-\frac{1}{it}\Im k_{a,\,\rho}(i))}
{1+(\frac{1}{it})^{2}} \\ = \frac{k_{a,\,\rho}(\frac{1}{it})-\Re
k_{a,\,\rho}(i)-\frac{1}{it}\Im k_{a,\,\rho}(i)}{it+(\frac{1}{it})}
\,\,. \ \ \ \ \ \ \ \ \
\end{multline*}
By letting $t= \frac{1}{w}>0$ we get
\begin{equation*}
\begin{split}
\int_{0}^{\infty}\hat{\rho}(\frac{s}{w})e^{-s}ds%
&= \frac{k_{a,\,\rho}(-iw)-\Re k_{a,\,\rho}(i)+iw \Im k_{a,\,\rho}(i)}{\frac{i}{w}-iw}\\
&= \frac{iwk_{a,\,\rho}(-iw)-iw\Re k_{a,\,\rho\mu}(i)-w^{2}\Im
k_{a,\,\rho}(i)}{w^{2}-1}
\end{split}
\end{equation*}
which, after substituting  $\frac{s}{w}=r>0$, is as follows
\[ \int_{0}^{\infty}\hat{\rho}(r)e^{-wr}dr=
\frac{iwk_{a,\, \rho}(-iw)-iw\Re k_{a, \,\rho}(i)-w^{2}\Im
k_{a,\,\rho}(i)}{w(w^{2}-1)}\,,
\]
and thus giving the formula in Theorem 1. Finally, inverting  the
Laplace transform of $\hat{\rho}$ and then inverting the Fourier
transform, we get uniquely the measure $\rho$ from values
$k_{a,\,\rho} (it), t\neq 0$. This completes the proof.

\emph{Proof of Corollary 1.} Simply note that
\begin{multline*}
\mathfrak{L}[\frac{1}{2}(e^{x}-e^{-x}); w]=\mathfrak{L}[\sinh x; w]
= \frac{1}{w^2-1},  \\ \mathfrak{L}[\frac{1}{2}(e^{x}+e^{-x});
w]=\mathfrak{L}[\cosh x; w] = \frac{w}{w^2-1}, \ \ w>1, \ \ \ \ \
\end{multline*}
which taken together with Theorem 1 give the proof.

\emph{Proof of Proposition 1.} Using the definitions (3) and (4) we
have
\begin{equation*}
\mathfrak{L}[\hat{\rho}; w]= \int_{\Rset}
\int_0^{\infty}\,e^{-r(w-ix)}\,
dr\,\rho(dx)=\int_{\Rset}\frac{1}{w-ix}\rho(dx)=-i\,g_{\rho}(-iw),
\end{equation*}
which completes the proof.

Here is an auxiliary lemma where the part (a) is a very standard
fact, recalled for completeness. This lemma simplifies the
arguments in the proof of the Example.
\begin{lem}
(a) If $P(z):=\prod_{j=1}^m\,(z-b_{j}), \ z\in\Cset$, for  some
complex numbers $b_j,\, j=1,2,..., m$, and $P'(z)$ is its derivative
then
\[
 \frac{P'(z)}{P(z)}\ =\sum_{j=1}^m\frac{1}{z-b_j}; \ \  \ \
 \frac{P''(z)}{P(z)}=(\sum_{j=1}^m \frac{1}{z-b_j}\,)^2-\sum_{j=1}^m\frac{1}{(z-b_j)^2};
\]
(b) If the $b_j$'s are distinct complex numbers and
$\xi_1,...,\xi_{m-1}$ denote the zeros of the equation $P'(z)=0$
then $\xi_j$ are different from $b_1,b_2,...,b_m$. Furthermore,
\begin{equation}
W_m(z):=\frac{z\,P'(z)-\,m\,P(z)}{P'(z)}=\frac{b_1+b_2+...+b_m}{m}+\sum_{j=1}^{m-1}\,\frac{\alpha_j}{z
-\xi_j}\,
\end{equation}
where
\[
\alpha_k:=  - m \frac{P(\xi_k)}{P''(\xi_k)}  =m\,\big[
\sum_{j=1}^m\frac{1}{(\xi_k-b_j)^2}- (\sum_{j=1}^m
\frac{1}{\xi_k-b_j})^2\big]^{-1}\,,
\]
for $k=1,2,...,m-1$.

(c) If the $b_j$'s are distinct real numbers, for $j=1,2,..., m$, then
$\alpha_k>0$, for $k=1,2,...,m-1$.
\end{lem}

\begin{proof}
(a)\ \ Since $P'(z)=\sum_{j=1}^{m}\prod_{k\neq j, k=1}^m\,(z-b_k)$
we get the first part of (a). Differentiating both sides of the
identity $P'(z)=P(z)\,\sum_{j=1}^m\frac{1}{z-b_j}$ we get the second
part of (a).

\noindent (b) \ \  Assume that $P$ and $P'$ have a common root.
Without loss of generality, lets say that $\xi_1=b_1$. Then
$P'(b_1)=\prod_{k=2}^m\,(b_1-b_k)=0$, which contradicts the
assumption that all $b_j$ are distinct.

Suppose that $\xi_1$ and its complex conjugate $\bar{\xi_1}$ are two
complex roots of $P'(z)=0$. Then from (a) we have
\[
P'(\xi_1)=P(\xi_1)\,\sum_{j=1}^{m-1}\frac{1}{\xi_1-b_j}=0=P(\bar{\xi_1})\,\sum_{j=1}^{m-1}\frac{1}{\bar{\xi_1}-b_j}\,.
\]
Since $P(\xi_1)\neq 0$ and $P(\bar{\xi_1})\neq 0$, we have
\[
\sum_{j=1}^{m-1}
\,[\frac{1}{\bar{\xi_1}-b_j}-\frac{1}{\xi_1-b_j}]=i\, 2 (\Im
\xi_1)\, \sum_{j=1}^m\frac{1}{|\xi_1-b_j|^2}=0 ,
\]
and hence $\Im\xi_1=0=\Im\xi_2=\Im\xi_3=...=\Im\xi_{m-1}$, that is,
all roots of $P'(z)=0$ are real.

\noindent Let us note that

 $P(z)=\prod_{k=1}^m(z-b_k)=z^m
+(-b_1-b_2-...-b_m)z^{m-1}+ Q_{m-2}(z)$,

\noindent for some polynomial $Q_{m-2}$ of degree m-2. Then
$z\,P'(z)-mP(z)= (b_1+...+b_m)z^{m-1}+\tilde{Q}_{m-2}(z)$ is a
polynomial of degree m-1, (for another polynomial of degree $m-2$).
Consequently, $W_m(z)$, given by (16), is a rational function (a
ratio of two polynomials of degree m-1). Since $\xi_1,...,\xi_{m-1}$
are zeros of $P'(z)=0$ , i.e., simple poles of $W_m(z)$, then
invoking the theorem on the decomposition of rational function into
a sum of simple fractions
\begin{multline}
W_m(z)=
z-\frac{m\,P(z)}{P'(z)}=\frac{(b_1+...+b_m)z^{m-1}+\tilde{Q}_{m-2}(z)}
{mz^{m-1}+(m-1)(-b_1-b_2-...-b_m)z^{m-2}+Q'_{m-2}(z)} \\ =
\frac{b_1+b_2+...+b_m}{m}+\sum_{j=1}^{m-1}\,\frac{\alpha_j}{z
-\xi_j} \ . \ \ \ \ \ \
\end{multline}
Putting $\bar{b}:= (b_1+...+b_m)/m$ and multiplying both sides by
$z-\xi_k$, we obtain
\[
\alpha_k+(z-\xi_k)\sum_{j\neq k , j =1}^{m-1}\,\frac{\alpha_j}{z
-\xi_j}=(z-\xi_k)(z-\bar{b}) -
m\,P(z)\Big(\frac{P'(z)-P'(\xi_k)}{z-\xi_k}\Big)^{-1}\,,
\]
and then letting $z\to\xi_k$ we explicitly get that
\[
\alpha_k:= - m \frac{P(\xi_k)}{P''(\xi_k)}= m\,\big[
\sum_{j=1}^m\frac{1}{(\xi_k-b_j)^2}- (\sum_{j=1}^m
\frac{1}{\xi_k-b_j})^2\big]^{-1}.
\]
  (c) Since $P(x)$ is a polynomial of m-th degree for $x\in\Rset$ and
$P(b_k)=P(b_{k+1})=0$ (for $b_j \in \Rset$) then, by the Mean Value
Theorem, there exists exactly one $\xi_j$ (in that interval) such
that $P'(\xi_k)=0$. If $P(\xi_k)>0$ then $P$ must be concave on that
interval and therefore $P''(\xi_k)<0$. Consequently, $\alpha_j>0$.
In the opposite case we have convex function that also leads to
the positivity of the $\alpha_k$ parameter. This completes the proof of
Lemma 1.
\end{proof}
\emph{Proof of the Example.} From Lemma 1 we have that the measure
$\rho_{\textbf{b}}$ is finite and positive. Furthermore, for
$a_{\textbf{b}}$ given by (6), using (16) (in Lemma 1) we get
\begin{multline*}
\int_{\Rset}\frac{1+zx}{z-x}d\rho_{\textbf{b}}(x)=
\sum_{j=1}^{m-1}\,\frac{1+\xi_j^2+z\xi_{j}-\xi_{j}^{2}}{z-\xi_{j}}\,\frac{\alpha_j}{1+\xi_j^2}=
\sum_{j=1}^{m-1} \frac{\alpha_j}{z-\xi_j}+\sum_{j=1}^{m-1}
\frac{\alpha_j\xi_j}{1+\xi^2_j}  \\
= W_m(z) - a_{\textbf{b}}=
z-\frac{m\,P_{\textbf{b}}(z)}{P'_{\textbf{b}}(z)}-
a_{\textbf{b}}=z-\frac{1}{G_{\mu_{\textbf{b}}(z)}}=E_{\mu_{\textbf{b}}}(z)-
a_{\textbf{b}}.
\end{multline*}
Substituting $it$ for $z$ in the above expression, we get equality (5) in the
Example.

\emph{Proof of Corollary 2.} Using (2) we obtain the expression (8)
for $-g_{\mu}(i)$. From (14) and (7), $e_{\mu}(i)=a-i\rho(\Rset)$,
we then infer the equalities in (9). [Note that
$d_{\mu}(1-d_{\mu})\ge c_{\mu}^2$].

In view of the assumption, $k_{a,\,\rho}$ in Corollary 1 may be
replaced by $e_{\mu}$, which combined with (7) and (9) yields
\[
ik_{a,\,\rho}(-iw)-i \Re k_{a,\, \rho}(i)-w \Im k_{a,\,\rho}(i)= w -
\frac{i}{g_{\mu}(-iw)} - i\frac{c_{\mu}}{|z_{\mu}|^2} - w
(1-\frac{d_{\mu}}{|z_{\mu}|^2}).
\]
Consequently the required identity follows from Corollary 1.

\emph{Proof of Proposition 2.} From Theorem 2.1 in Chistyakov and
Goetze (2005), (cf. also Biane (1998)) for measures $\mu_1$ and
$\mu_2$ there exist uniquely determined probability measures
$\nu_1$, $\nu_2$ and $\mu$ such that
\[
z=F_{\nu_1}(z)+ F_{\nu_2}(z)-F_{\mu_1}(F_{\nu_1}(z)) \ \ \ \
\mbox{and} \ \
F_{\mu_1}(F_{\nu_1}(z))=F_{\mu_2}(F_{\nu_2}(z))=F_{\mu}(z),
\]
where $\mu=\mu_1\boxplus\mu_2$ (Voiculescu convolution). Hence
\[
E_{\nu_1\uplus \nu_2}(z)= E_{\nu_1}(z) + E_{\nu_2}(z)=
z-F_{\nu_1}(z)+ z-
F_{\nu_2}(z)=z-F_{\mu_1\boxplus\mu_2}=E_{\mu_1\boxplus\mu_2}(z).
\]
From the uniqueness of the self-energy functional we get
$\mu_1\boxplus\mu_2={\nu_1\uplus \nu_2}$, which completes the proof.

\emph{Proof of Corollary 3.} From Corollary 2.2 in Chistyakov and
Goetze (2005), for measures $\mu_1,...,\mu_n$ there exist uniquely
determined probability measures $\nu_1,...,\nu_n$ and $\mu$ such
that
\begin{multline*}
z=F_{\nu_1}(z)+ ... + F_{\nu_n}(z)- (n-1)F_{\mu_1}(F_{\nu_1}(z)) \ \
\ \ \ \ \ \ \\ \mbox{and} \ \
F_{\mu_1}(F_{\nu_1}(z))=...=F_{\mu_n}(F_{\nu_n}(z))=F_{\mu}(z),  \ \
\ \ \ \ \
\end{multline*}
where $\mu=\mu_1\boxplus ... \boxplus\mu_n$ (the Voiculescu
convolution). Thus
\begin{multline*}
E_{\nu_1\uplus . . . \uplus\nu_2}(z)= E_{\nu_1}(z)+...+
E_{\nu_n}(z)= z-F_{\nu_1}(z)+ ...+ z- F_{\nu_n}(z) \\ =
(n-1)(z-F_{\mu_1}(F_{\nu_1}(z))=
(n-1)(z-F_{\mu_1\boxplus...\boxplus\mu_n}(z)
=(n-1)E_{\mu_1\boxplus...\boxplus \mu_n},
\end{multline*}
which completes the proof.

\medskip
\medskip
\textbf{Acknowledgements.} We would like to thank anonymous Reviewer
who called our attention to the paper by R. Lenczewski (2007), which
contains much more then our Proposition 2 and Corollary 3. In
particular, he identified new s-free independence, new convolution
and found the corresponding Hilbert space decomposition. The second
named co-author thanks R. Lenczewski for a fruitful personal
discussion.

\newpage
\begin{center}
\textbf{REFERENCES}
\end{center}

\noindent [1] N. I. Akhiezer (1965) \emph{The classical moment
problem}, Oliver \& Boyd, Edinburgh and London

\noindent [2] H. Bercovici, D. Voiculescu (1993) Free convolution of
measures with unbounded support, \emph{Indiana Univ. Math. J.}
vol.42, No.3, pp. 733-773.

\noindent [3] Ph. Biane (1998). Processes with free increments,
\emph{Math. Z.}, vol. 227, pp. 143-174.

\noindent [3] L. Bondesson (1992). \emph{Generalized gamma
convolutions and related} \emph{classes of distributions and
densities} Springer-Verlag, New York.

 \noindent[4] G. P. Chistyakov,
F. $G\ddot{o}tze$ (2005). \emph{The arithmetic of distributions in
free probability theory}, arXiv:math.OA/0508245

\noindent [5] I. S. Gradshteyn and I.M. Ryzhik (1994), \emph{Tables
of integrals, series and products}, Academic Press, San Diego; 5th
Edition.

\noindent [6] Z. J. Jurek (2006). Cauchy transforms of measures
viewed as some functionals of Fourier transforms, \emph{Probab.
Math. Stat.}, vol. 26 Fasc. 1, pp 187-200.

\noindent [6] Z. J. Jurek (2007). Random integral representations
for free-infinitely divisible and tempered stable distributions,
\emph{Stat.$\&$ Probab. Letters}, vol. 77, pp. 417-425.

\noindent [7] S. Lang (1975), \emph{$SL_2(\Rset)$}, Addison-Wesley,
Reading Massachusetts.

\noindent [8] R. Lenczewski (2007). Decompositions of the free
additive convolution, \emph{J. Funct. Analysis} vol. 246, pp.
330-365.

\noindent [9] R. Speicher and R. Woroudi (1997). \emph{Boolean
convolution} Fields Institute Communications, vol.12, pp 267-279

\medskip
\noindent Institute of Mathematics, University of Wroc\l aw,
Pl.Grunwaldzki 2/4, 50-384 Wroc\l aw, Poland \\
e-mail: zjjurek@math.uni.wroc.pl
 \ \ www.math.uni.wroc.pl/$^{\sim}$zjjurek

\end{document}